\documentclass[11pt,a4paper]{article}
\usepackage{latexsym,amssymb,amsfonts,amsmath,amsthm,nccmath,float,enumitem,setspace,bm}
\usepackage[hidelinks]{hyperref}
\usepackage[dvipsnames]{xcolor}
\usepackage[margin=2cm]{geometry}
\usepackage{tikz}
\usetikzlibrary{arrows}
\usetikzlibrary{calc}
\usepackage{color}

\newtheorem{theorem}{Theorem}

\newtheorem{conjecture}[theorem]{Conjecture}
\theoremstyle{definition}

\date{}
\begin{document}

\title{Counterexamples to a conjecture on matching Kneser graphs}

\author{Moharram N. Iradmusa\thanks{Faculty of Mathematical Sciences, Shahid Beheshti University, Tehran, Iran ({\tt m\_iradmusa@sbu.ac.ir})}}
\maketitle

\begin{abstract}
Let $G$ be a graph and $r\in\mathbb{N}$. The matching Kneser graph $\textsf{KG}(G, rK_2)$ is a graph whose vertex set is the set of $r$-matchings in $G$ and two vertices are adjacent if their corresponding matchings are edge-disjoint. In [Alishahi, M. and Hajiabolhassan, H., On the Chromatic Number of Matching Kneser Graphs, Combin. Probab. and Comput. 29 (2020), no. 1, 1--21.] it was conjectured that for any connected graph $G$ and positive integer $r\geq 2$, the chromatic number of $\textsf{KG}(G, rK_2)$ is equal to $|E(G)|-\textsf{ex}(G,rK_2)$, where $\textsf{ex}(G,rK_2)$ denotes the largest number of edges in $G$ avoiding a matching of size $r$. In this note, we show that the conjecture is not true for snarks.
\end{abstract}

%\setstretch{1.1}

%\section{A conjecture on matching Kneser graphs}
Given two positive integers $n$ and $r$, the Kneser graph $\textsf{KG}(n,r)$, is the graph whose vertices represent the $r$-subsets of $\{1,\ldots,n\}$, and where two vertices are adjacent if and only if they correspond to disjoint subsets. This graph introduced by Lov\'{a}sz in 1978 \cite{lovasz}. The matching Kneser graph was introduced as a generalization of Kneser graph \cite{alishahi}. Precisely, $\textsf{KG}(nK_2,rK_2)\cong\textsf{KG}(n,r)$. Also, by a simple greedy coloring, it was proved that $\chi(\textsf{KG}(G,rK_2))\leq |E(G)|-\textsf{ex}(G,rK_2)$ for any graph $G$ and positive integer $r$. This upper bound  is tight when $G$ satisfies some certain properties as is claimed in the main result of \cite{alishahi} and is shown throughout the paper.
Finally the following conjecture was introduced according to the results of the paper.
\begin{conjecture}
For any connected graph $G$ and positive integer $r\geq 2$,
\[\chi(\textsf{KG}(G,rK_2))=|E(G)|-\textsf{ex}(G,rK_2).\] 
\end{conjecture}
%\section{A counterexample}
We show that the conjecture is false when $G$ is a snark. A snark is a simple, connected, bridgeless cubic graph with chromatic index equal to four \cite{gardner}. Let $G$ be a snark of order $2r$. Petersen's Theorem from 1891 \cite{petersen}, states that every bridgeless cubic graph has a perfect matching. Further, Sch\"{o}nberger in 1934 has proved that every bridgeless cubic graph has a perfect matching not containing two arbitrarily prescribed edges \cite{schonberger}. Therefore, $\textsf{ex}(G,rK_2)\leq |E(G)|-3$. In addition, by removing three edges incident with an arbitrary vertex, the resulting graph has no matching of size $r$. Therefore $\textsf{ex}(G,rK_2)=|E(G)|-3$. In addition, since $\chi'(G)=4$, the intersection of any two perfect matchings of $G$ is nonempty. Therefore, $\textsf{KG}(G,rK_2)$ has no edge and so $\chi(\textsf{KG}(G,rK_2))=1\neq |E(G)|-\textsf{ex}(G,rK_2)=3$.

%\noindent\textbf{Acknowledgments.} The first author was supported by Australian Research Council grants DP150100506 and FT160100048.

\end{document}